\documentclass[12pt]{article}
\usepackage{cite}
\usepackage{amssymb}
\usepackage{amsmath}
\usepackage{graphicx}
\newtheorem{thm}{Theorem}[section]
\newtheorem{cor}[thm]{Corollary}

\usepackage{graphics}
\numberwithin{equation}{section}

\newcommand{\eh}{\hfill}\newlength{\sperr}

\def\ve{\varepsilon}

\newtheorem{Pa}{Paper}[section]
\newtheorem{Tm}[Pa]{{\bf Theorem}}

\newtheorem{Rk}[Pa]{{\bf Remark}}

\newtheorem{Dn}[Pa]{{\bf Definition}}
\newtheorem{Pn}[Pa]{{\bf Proposition}}
\numberwithin{equation}{section}

\title{Levy processes with summable   Levy measures, long time behavior}
\author{Lev Sakhnovich}
\date{}
\begin{document}
\maketitle

\thanks{99 Cove ave., Milford, CT, 06461, USA \\
 E-mail: lsakhnovich@gmail.com}\\

\begin{abstract} In our previous  paper \cite{Sakh7} we have proved  that a  representation
of the infinitesimal generators $L$ for Levy processes $X_t$ can be written down in a
convolution type form. For the case of non-summable Levy measures we constructed
the quasi-potential operators $B$ and investigated the long time behavior of $X_t$.
In the present paper we consider  Levy processes $X_t$ with summable
Levy measures. In this case the form of the quasi-potential operators $B$  essentially differs from the form
in the case of non-summable Levy measures. We use this new form in order to study the long time behavior of $X_t$
for the case of summable
Levy measures.
  \end{abstract}

\section{Introduction}\label{sec1}
%%%%%%%%

 Let us introduce  the notion of the Levy processes.

\begin{Dn}\label{Definition 1.1.}
A stochastic process $\{X_{t}:t{\geq}0\}$ is called Levy
process, if the following conditions are fulfilled:\\
1. Almost surely $X_{0}=0$, i.e. $P(X_{0}=0)=1$.\\
(One says that an event happens almost surely (a.s.) if it happens with probability one.)\\
2. For any $0{\leq}t_{1}<t_{2}...<t_{n}<\infty$ the random variables \\
 $X_{t_{2}}-X_{t_{1}}, X_{t_{3}}-X_{t_{4}},..., X_{t_{n}}-X_{t_{n-1}}$\\ are independent (independent increments).\\
( To call the increments of the process $X_{t}$ \emph{independent} means that
increments $X_{t_{2}}-X_{t_{1}}, X_{t_{3}}-X_{t_{4}},..., X_{t_{n}}-X_{t_{n-1}}$ are mutually  (not just
pairwise) independent.)\\
3. For any $s<t$ the distributions of   $X_{t}-X_{s}$ and $X_{t-s}$ are equal
(stationary increments).\\
4. Process $X_{t}$ is almost surely right continuous with left limits.\\
\end{Dn}
Then Levy-Khinchine formula gives (see \cite{Bert}, \cite{Sato})
 \begin{equation}
\mu(z,t)=E\{\mathrm{exp}[izX_{t}]\}=
\mathrm{exp}[-t\lambda(z)],\quad t{\geq}0,  \label{1.1} \end{equation}
where
\begin{equation}
\lambda(z)=\frac{1}{2}Az^{2}-i{\gamma}z-\int_{-\infty}^{\infty}(e^{ixz}-1-ixz1_{|x|<1})\nu(dx).
 \label{1.2} \end{equation} Here $A{\geq}0,\quad \gamma=\overline{\gamma},\quad z=\overline{z}$ and
 $\nu(dx)$ is a measure on the axis $(-\infty,\infty)$
  satisfying the conditions
\begin{equation}
\int_{-\infty}^{\infty}\frac{x^{2}}{1+x^{2}}\nu(dx)<\infty.
 \label{1.3} \end{equation}
The  Levy-Khinchine formula is determined by the Levy-Khinchine triplet
(A,$\gamma, \nu(dx)$).\\
By $P_{t}(x_{0},\Delta)$ we denote the probability
$P(X_{t}{\in}\Delta)$ when $P(X_{0}=x_{0})=1$ and $\Delta{\in}R$.
The transition operator $P_{t}$ is defined by the formula
\begin{equation}
P_{t}f(x)=\int_{-\infty}^{\infty}P_{t}(x,dy)f(y). \label{1.4} \end{equation} Let
$C_{0}$ be the Banach space of continuous functions $f(x)$ ,
satisfying the condition $\mathrm{lim}f(x)=0,\quad |x|{\to}\infty$
with the norm $||f||=\mathrm{sup}_{x}|f(x)|$. We denote by
$C_{0}^{n}$ the set of $f(x){\in}C_{0}$ such that
$f^{(k)}(x){\in}C_{0},\quad (1{\leq}k{\leq}n).$ It is known that
\cite{Sato}
\begin{equation} P_{t}f{\in}C_{0}, \label{1.5} \end{equation}
if $f(x){\in}C_{0}^{2}.$\\
Now we formulate the following important result (see \cite{Sato}) .
\begin{Tm}\label{Theorem 1.2.}
{The family of the operators $P_{t}\quad
(t{\geq}0)$ defined by the Levy process $X_{t}$ is a strongly
continuous  semigroup on $C_{0}$ with the norm $||P_{t}||=1$. Let
$L$ be its infinitesimal generator. Then}
\begin{equation} Lf=\frac{1}{2}A\frac {d^{2}f}{dx^{2}}+\gamma
\frac{df}{dx}+\int_{-\infty}^{\infty}(f(x+y)-f(x)-y\frac{df}{dx}1_{|y|<1})\nu(dy), \label{1.6} \end{equation}
{where} $f{\in}C_{0}^{2}$.
\end{Tm}
Slightly changing the well-known classification (see\cite{Sato}) we introduce the following definition:\\
\begin{Dn}\label{Definition 1.3} We say that a Levy process $X_{t}$ generated by $(A,\nu,\gamma)$ has  type
$I$ if
\begin{equation}A=0 \,and\,\int_{-\infty}^{\infty}\nu(dx)<\infty,\label{1.7}\end{equation}
and $X_{t}$ has  type $II$ if
\begin{equation}A{\ne}0 \quad or\quad \int_{-\infty}^{\infty}\nu(dx)=\infty.\label{1.8}\end{equation}
\end{Dn}
\begin{Rk}\label{Remark 1.4}The introduced type $I$ coincides with the type $A$ in the
 usual classification.The introduced type $II$ coincides with the union  of the types $B$ and $C$ in the usual classification.\end{Rk}
The properties of these two types of the Levy processes are quite  different.
The paper \cite{Sakh7} was  dedicated to type $II$.
\begin{Rk}\label{Remark 1.5}Some results of the paper \cite{Sakh7} are proved for all Levy processes.
 We shall use these results.\end{Rk}
 In the present paper
 we shall consider the type $I$. Without loss of generality we assume that
 \begin{equation}\gamma-\int_{|y|<1}y\nu(dy)=0.\label{1.9}\end{equation}
  According to \eqref{1.6} the generator $L$ of the corresponding process $X_t$ can be represented in the form (see \cite{Sakh6},p.13)
 \begin{equation}
Lf=\int_{-\infty}^{\infty}[f(x+y)-f(x)]\nu(dy).\label{1.10}\end{equation}
\begin{Rk}\label{Remark 1.6}If condition \eqref{1.9} is valid, then type $I$ coincides with the class of the
compound
Poisson processes.\end{Rk}
As in the case of type $II$ we use the convolution representation of the generator $L$.
To do it
 we introduce the functions
\begin{equation}\mu_{-}(x)=\int_{-\infty}^{x}\nu(dx),\,x<0, \label{1.11}\end{equation}
\begin{equation}\mu_{+}(x)=-\int_{x}^{\infty}\nu(dx),\,x>0,\label{1.12} \end{equation}
where the functions $\mu_{-}(x)$ and $\mu_{+}(x)$ are  monotonically increasing and right continuous on the  half-axis   $(-\infty,0]$ and $[0,\infty)$ respectively. We note that
\begin{equation}\mu_{+}(x){\to}0,\,x{\to}+\infty;\,
\mu_{-}(x){\to}0,\,x{\to}-\infty,\label{1.13}\end{equation}
\begin{equation}\mu_{-}(x){\geq}0,\,x<0;\,\mu_{+}(x){\leq}0,\,x>0.\label{1.14}\end{equation}
Now we define  the functions
\begin{equation}k_{-}(x)=\int_{-a}^{x}\mu_{-}(t)dt,\,-\infty{\leq}x<0,\,a>0,\label{1.15}
\end{equation}
\begin{equation}k_{+}(x)=-\int_{x}^{a}\mu_{+}(t)dt,\,0<x{\leq}+\infty.\label{1.16}
\end{equation}
\begin{Pn}\label{Proposition 1.7}(see \cite{Sakh6},p.13). Let the Levy process $X_t$ belong to the type $I$ and let the condition \eqref{1.9}
be fulfilled. Then  formula \eqref{1.10} can be written in the following
convolution form
\begin{equation}Lf=\frac{d}{dx}S\frac{d}{dx}f,\label{1.17}\end{equation}
where
\begin{equation}Sf=\int_{-\infty}^{\infty}k(y-x)f(y)dy,\label{1.18}\end{equation}
\begin{equation}k(x)=k_{+}(x),\,if\,x>0;\,k(x)=k_{-}(x),\,if\,x<0.\label{1.19}\end{equation}
\end{Pn}
Using  \eqref{1.15}-\eqref{1.18} we obtain the assertion.\\
\begin{Pn}\label{Proposition 1.8} Let conditions of Proposition 1.7\, be fulfilled. Then relation (1.17)
takes the form
\begin{equation}Lf=-{\Omega}f+\int_{-\infty}^{\infty}f(y)d_{y}\nu(y-x),\,
\Omega=\int_{-\infty}^{\infty}d\nu(y).\end{equation}\label{1.20}\end{Pn}
\emph{Proof.} It follows from \eqref{1.15}-\eqref{1.18} that
\begin{equation}Lf=-\int_{-\infty}^{x}\mu_{-}(y-x)f^{\prime}(y)dy
-\int_{x}^{\infty}\mu_{+}(y-x)f^{\prime}(y)dy.\label{1.21}\end{equation}
Integrating by parts \eqref{1.21} we have
\begin{equation}Lf=-[\mu_{-}(0)-\mu_{+}(0)]f+\int_{-\infty}^{\infty}f(y)d_{y}\nu(y-x).
\label{1.22}\end{equation}The proposition is proved.\\
In formulas \eqref{1.20} and \eqref{1.22} we use the equality $d\nu(x)=\nu(dx)$.\\
\begin{Dn}\label{Definition 1.9}We say that Levy process $X_t$ belongs to type $I_c$ if the corresponding
Levy measure $\nu(y)$ is summable and continuous.\\
We say that Levy process $X_t$ belongs to type $I_d$ if the corresponding
Levy measure $\nu(y)$ is summable and discrete.\end{Dn}
It is easy to see that the following assertion is true.
\begin{Pn}\label{Proposition 1.10}If Levy process $X_t$ belongs to the type $I$, then
 $X_t$ can be represented in the form
\begin{equation}X_t=X^{(1)}_t+X^{(2)}_t, \label{1.23}\end{equation}
where $X^{(1)}_t{\in}I_c$ and $X^{(2)}_t{\in}I_d$.\end{Pn}
 The main part of the present paper is dedicated to investigating  the Levy processes from the type $I_c$.\\
We denote by $p(t,\Delta)$ the probability  that a sample of the process $X_t{\in}I_c$ remains inside
the domain $\Delta$ for $0{\leq}\tau{\leq}t$ (ruin problem). With the help of representation \eqref{1.17} we find a new
formula for $p(t,\Delta)$.
This formula allow us to obtain the long time behavior of $p(t,\Delta)$. Namely, we have proved the following asymptotic formula
\begin{equation}
 p(t,\Delta)=e^{-t/\lambda_{1}}[c_{1}+o(1)],\,c_1>0,\,\lambda_{1}>0,\,\quad t{\to}+\infty.\label{1.24}
   \end{equation}
Let $T_{\Delta}$  be the time  during which $X_{t}$ remains in the domain $\Delta$
before it leaves the domain $\Delta$ for the first time. It is easy to see that
\begin{equation}p(t,\Delta)=P(T_{\Delta}>t).\label{1.25}\end{equation}
An essential role in our theory plays the operator
\begin{equation}L_{\Delta}f=-{\Omega}f+\int_{\Delta}f(y)d_{y}\nu(y-x),
\label{1.26}\end{equation}which is generated by the operator $L$ (see \eqref{1.20}).
We note that $\lambda_{1}$ in formula \eqref{1.23} is the greatest eigenvalue of $-L_{\Delta}^{-1}$.
\begin{Dn}\label{Definition 1.11}The measure $\nu(y)$ is unimodal with mode $0$ if $\nu(y)$ is concave when
$y<o$ and convex if $y>0$.\end{Dn}The unimodality and it properties were actively investigated (see\cite{Sato}). In the paper we found a new important property
of unimodal measure:
\begin{Pn}\label{Proposition 1.12}If Levy measure $\nu(y)$ is continuous, summable and unimodal with mode {0} then the operator $L_{\Delta}^{-1}$ has the form
\begin{equation}L_{\Delta}^{-1}=-\frac{1}{\Omega}(I+T_{1}),\label{1.27}\end{equation}
where the operator $T_{1}$ is compact in the space of continuous functions.\end{Pn}
In the last part of the paper we investigate the operator $L_{\Delta}$ when $X_t{\in}I_d$.

%%%
\section{Quasi-potential}
\label{section 2}
1.By domain $\Delta$ we denote  the set of segments $[a_{k},b_{k}]$ , where\\
 $a_{1}<b_{1}<a_{2}<b_{2}<...<a_{n}<b_{n},\quad 1{\leq}k{\leq}n.$\\
We denote by $D_{\Delta}$ the space of the continuous functions $g(x)$ on the domain $\Delta$.
The norm in $D_{\Delta}$ is defined by the relation $||f||=sup_{x{\in}\Delta}|f(x)|.$
The space $D_{\Delta}^{0}$ is defined by the relations:\\
$f(x){\in}D_{\Delta}$ and $f(a_k)=f(b_k)=0,\quad (1{\leq}k{\leq}n).$|\\
 We introduce the operator $P_{\Delta}$ by relation
$P_{\Delta}f(x)=f(x)$ if $x{\in}\Delta$ and $P_{\Delta}f(x)=0$ if $x{\notin}\Delta.$\\
\begin{Dn}\label{Definition 2.1}
The operator\begin{equation}
L_{\Delta}=P_{\Delta}LP_{\Delta}=\frac{d}{dx}S_{\Delta}\frac{d}{dx},
\,where\,S_{\Delta}=P_{\Delta}SP_{\Delta},  \label{2.1} \end{equation}
is called a  truncated generator. (We use here
the equality $P_{\Delta}\frac{d}{dx}=\frac{d}{dx}{P_{\Delta}},\,x{\in}\Delta).$
\end{Dn}
\begin{Dn}\label{Definition 2.2.}
The operator $B$ with the definition domain
 $D_\Delta$  is called a quasi-potential if the
following relation
\begin{equation} -BL_{\Delta}g=g,\quad g{\in}D^0_{\Delta} \label{2.2} \end{equation}
is true.\end{Dn}
According to Proposition 1.8 the operator $L_{\Delta}$ has the form
\begin{equation}L_{\Delta}f=-{\Omega}f+\int_{\Delta}f(y)d_{y}\nu(y-x),\,f(x){\in}D_{\Delta}^{0},
\label{2.3}\end{equation}We introduce the operator
\begin{equation}Tf=\int_{\Delta}f(y)d_{y}\nu(y-x),\,f(x){\in}D_{\Delta}^{0}.\label{2.4}
\end{equation}Inequality \eqref{2.4} implies the statement.
\begin{Tm}\label{Theorem 2.3}Let the Levy process $X_t$ belong to the type $I$
and let condition \eqref{1.9} be fulfilled. Then the operator $T$ acts from $D_{\Delta}$
into $D_{\Delta}$ and
\begin{equation}||T||=sup{\int_{\Delta}d_{y}\nu(y-x)}{\leq}\Omega,\quad x{\in}\Delta
\label{2.5}\end{equation}\end{Tm}
Further we suppose in addition that
\begin{equation}||T||<\Omega.\label{2.6}\end{equation}
\begin{Rk}\label{Remark 2.4} Let the support of the   Levy measure $\nu$ is unbounded.Then in  view of \eqref{2.5}
the inequality \eqref{2.6} holds.\end{Rk}
 We consider separately the case when the Levy process $X_t$ belongs to the type $I_c$.
\begin{Tm}\label{Theorem 2.5}Let $X_t$ belong to the type $I_c$ and let the condition \eqref{2.6} be true.\\
Then the operator $L^{-1}_{\Delta}$ exists and has the form
\begin{equation}L^{-1}_{\Delta}=-B=-\frac{1}{\Omega}(I+T_{1}),\label{2.7}\end{equation}
where the operator $T_1$ acts in the space $D_{\Delta}$ is bounded and is defined by the formula
\begin{equation}T_{1}f=\int_{\Delta}f(y)d_{y}\Phi(x,y).\label{2.8}\end{equation}
The function $\Phi(x,y)$ is continuous with respect to $x$ and $y$, monotonically
increasing with respect to $y$.\end{Tm}
\section{Quasi-potential, compactness}
1. In this section  we consider the following problem:\\
Under which conditions the operators $T$ and $T_1$ are compact in the space $D_{\Delta}$?\\
We remind that the operators $T$ and $T_1$ are defined by \eqref{2.4} and
\eqref{2.8} respectively. We need the following definitions
\begin{Dn}\label{Definition 3.1} The total variation of a complex-valued function g, defined on $\Delta$ is the quantity\\
$$V_{\Delta}(g)=\sup_{P}\sum_{i=0}^{n_P-1}|g(x_{i+1})-g(x_{i})|,$$\\
where the supremum is taken over the set  of all partitions $P=(x_0,x_1,...,x_{n_P})$
of the $\Delta$.\end{Dn}
\begin{Dn}\label{Definition 3.2}. A complex-valued function g on the  $\Delta$ is said to be of bounded variation (BV function) on the $\Delta$ if its total variation is finite.
\end{Dn} By $D^*_{\Delta}$
we denote the conjugate space to $D_{\Delta}$. It is well-known that the space $D^*_{\Delta}$ consists from functions $g(x)$ with a bounded total variation $V_{\Delta}(g)$. The norm in $D^*_{\Delta}$ is defined by the relation $||g||=V_{\Delta}(g)$, the functional in $D_{\Delta}$ is defined by the relation
\begin{equation}(f,g)_{\Delta}=\int_{\Delta}f(x)d\overline{g(x)},\,f{\in}D_{\Delta},\,
g{\in}D^*_{\Delta}.\label{3.1}\end{equation} Hence, the conjugate operator
$B^*$ maps the space $D^*_{\Delta}$ into himself and has the form
\begin{equation}B^*g=\int_{\Delta}\Phi(y,x)dg(y).\label{3.2}\end{equation}

J.Radon \cite{RAD} proved the following theorem.
\begin{Tm}\label{Theorem 3.3} The operator $T_{1}$ defined by formula \eqref{2.8} is compact in the space $D_{\Delta}$
if and only if
\begin{equation}\lim_{x{\to}\xi}||\Phi(x,y)-\Phi(\xi,y)||_{V}=0,\,x,y,\xi{\in}\Delta. \end{equation}\label{3.3}\end{Tm}
Hence we have the assertion.
\begin{Pn}\label{Proposition 3.4}If measure $\nu(y)$ is summable and has continuous derivative,
then the corresponding operators $T$ and $T_1$ are compact in the space $D_{\Delta}$.
\end{Pn}
\emph{Proof.} Relation  \eqref{2.4} takes the form:
\begin{equation}Tf=\int_{\Delta}f(y)\nu^{\prime}(y-x)dy,\,x{\in}\Delta.\label{3.4}\end{equation}
The conditions of the Radon's theorem are fulfilled, i.e. the operator $T$ is compact.
In view of \eqref{2.4},\eqref{2.5} and \eqref{2.8} the operator $T_1$ is compact as well.
The Proposition is proved.\\
2.Let us  consider the important case, when the Levy measure is unimodal (see Definition 1.11).
We shall use the following convex and concave properties.
\begin{Pn}\label{Proposition 3.5} Let the points $x_p,\,(p=1,2.3.4)$ be such that
$x_1<x_2{\leq}x_3<x_4 $.\\
1. If a function $f(x)$ is convex then
\begin{equation}\frac{f(x_2)-f(x_1)}{x_2-x_1}{\geq}\frac{f(x_4)-f(x_3)}{x_4-x_3}.
\label{3.5}\end{equation}
2. If a function $f(x)$ is concave then
\begin{equation}\frac{f(x_2)-f(x_1)}{x_2-x_1}{\leq}\frac{f(x_4)-f(x_3)}{x_4-x_3}.
\label{3.6}\end{equation}\end{Pn}
\begin{cor}\label{corollary 3.1}The assertions 1. and 2. of Proposition 3.5 are true  if\\
$x_1<x_3{\leq}x_2<x_4 $ and $x_2-x_1=x_4-x_3 $.
\end{cor}
\emph{Proof.} We consider the convex case and    take such integer $n$ that $\ell/n<\ell_1$, where $\ell=x_2-x_1,$
$\ell_1=x_3-x_1$.It follows  from  (3.5) that
\begin{equation}\sum_{k=1}^{n}[f(x_1+k\ell/n)-f(x_1+(k-1)\ell/n)]{\geq}
\sum_{k=1}^{n}[f(x_3+k\ell/n)-f(x_3+(k-1)\ell/n)].\label{3.7}\end{equation}Hence in the convex case the corollary is proved. In the same way the corollary can be proved in the concave case.
\begin{Tm}\label{Theorem 3.6}If a measure $\nu(dx)$ on $R$ is unimodal with mode $0$ then the corresponding
operators $T$ and $T_1$ are compact in the space $D_{\Delta}$.\end{Tm}
\emph{Proof.}

 Let us consider the case when $\Delta=[c,d]$ and\\
$c=y_0<y_1<...<y_n=d.$ We introduce the variation
\begin{equation}V_{n}(x,\xi)=\sum_{k=1}^{n}|[\mu(y_k-x)-\mu(y_{k-1}-x)]-[\mu(y_k-\xi)-\mu(y_{k-1}-\xi)]|,
\label{3.8}\end{equation} where $c{\leq}\xi<x{\leq}d.$ Without loss of generality we assume that
\begin{equation}max|y_k-y_{k-1}|<x-\xi,\,1{\leq}k{\leq}n.\label{3.9}\end{equation}
We denote by $y_N$ such point that
\begin{equation}y_N-x{\geq}0, \, y_{N-1}-x{\leq}0.\label{3.10}\end{equation}We represent equality  \eqref{3.8} in the form
\begin{equation}V_{n}(x,\xi)=\sum_{k=1}^{N-2}|b_k|+\sum_{k=N-1}^{N}|b_k|+\sum_{N+1}^{n}|b_k|,
\label{3.11}\end{equation}
where
\begin{equation}b_k=[\mu(y_k-x)-\mu(y_{k-1}-x)]-[\mu(y_k-\xi)-\mu(y_{k-1}-\xi)].\label{3.12}
\end{equation}
 Proposition 3.6 implies that
\begin{equation}[\mu(y_k-x)-\mu(y_{k-1}-x)]-[\mu(y_k-\xi)-\mu(y_{k-1}-\xi)]{\geq}0,\,k{\geq}N,
\label{3.13}\end{equation}
\begin{equation}[\mu(y_k-x)-\mu(y_{k-1}-x)]-[\mu(y_k-\xi)-\mu(y_{k-1}-\xi)]{\leq}0,\, k{\leq}N-2.\label{3.14}\end{equation} It follows from \eqref{3.11}-\eqref{3.14} that
\begin{equation}V_{n}(x,\xi)=-\sum_{k=1}^{N-2}b_k+\sum_{k=N-1}^{N}|b_k|+\sum_{N+1}^{n}b_k.
\label{3.15}\end{equation}Hence we have
\begin{equation}V_{n}(x,\xi){\leq}D_0+2(D_{N-2}+D_{N-1}+D_{N})+D_n,\label{3.16}\end{equation}
where $D_k=|\mu(y_k-x)-\mu(y_k-\xi)|.$ The function $\mu(y)$ is continuous. Therefore
\begin{equation}supV_{n}(x,\xi){\to}0,\,x{\to}\xi.\label{3.17}\end{equation}From the last relation and Radon's  theorem follows that in the case ,when $\Delta=[c,d],$ the  assertion of the theorem is true.
Then the theorem is true in the case  of arbitrary domain $\Delta$.
\begin{Rk}\label{Remark 3.8}The Theorem 3.6 is true in the case when $\nu(dx)$ is n-modal $(1{\leq}n<\infty)$.\end{Rk}
\begin{Rk}\label{Remark 3.9}The unimodality of the Levy measure $\nu(dx)$ is closely  connected with
the  unimodality
of the probability distribution $F(x,t)$ of the corresponding Levy process ( see \cite{Sato},section 52).\end{Rk}
It is easy to obtain the following assertion.\\
\begin{Pn}\label{Proposition 3.10}Let condition \eqref{2.6} be fulfilled.Then the spectrum of the operator $B$
belongs to the right half plane.If in addition the operator $T$ is compact, then
the  eigenvalues $\lambda_j$ of $B$ are such that
\begin{equation}\lambda_j{\to}1/{\Omega},\, j{\to}\infty.\label{3.18}\end{equation}\end{Pn}
\section{The Probability of the Levy process (type $I_c$) remaining within the given domain}
\label{sec 4}
%%%%%
1.  We remind, that the definition of the Levy processes and the definition of the type  $I_c$ are given
in  section 1.

\textbf{Conditions 4.1} \label{Conditions 4.1}\emph{Further we assume that the following conditions are
fulfilled:\\
1.The Levy process $X_t$ belong to the type $I_c$. \\
2.The relation \eqref{1.9} is valid.}\\
We denote by $F_{0}(x,t)$  the distribution function of Levy process $X_t$, i.e.
   \begin{equation}F_{0}(x,t)=P(X_{t}{\leq}x).
 \label{4.1}\end{equation}We need the following statement (see \cite{Sato}, Remark 27.3)
\begin{Tm}\label{Theorem 4.1}Let Conditions 4.1 be fulfilled. The distribution function $F_{0}(x,t)$
is continuous with respect to $x$ if $x{\ne}0$. \end{Tm}
Let us investigate the behavior of  $F_{0}(x,t)$ in the point $x=0$.
\begin{Pn}\label{Proposition 4.2} Let Conditions 4.1 be fulfilled. Then the relation
\begin{equation}F_{0}(+0,t)-F_{0}(-0,t)=e^{-t\Omega},\,
\Omega=\int_{-\infty}^{\infty}\nu(dx) \label{4.2}\end{equation} is valid.
Here by definition we have
\begin{equation}F(0,t)=[F_{0}(+0,t)+F_{0}(-0,t)]/2.\label{4.3})\end{equation}\end{Pn}
\emph{Proof.} In our case equality \eqref{1.2} takes the form
\begin{equation}\lambda(z)=-\int_{-\infty}^{\infty}(e^{ixz}-1)\nu(dx)=\Omega-\omega(z),
\label{4.4}\end{equation}where
\begin{equation}\omega(z)=\int_{-\infty}^{\infty}e^{ixz}\nu(dx).\label{4.5}\end{equation}
Using the inverted Fourier-Stieltjes transform we have
\begin{equation}F_{0}(x,t)-F_{0}(0,t)=\frac{1}{2\pi}e^{-t\Omega}
\int_{-\infty}^{\infty}\frac{e^{-ixz}-1}{-iz}
\sum_{k=0}^{\infty}\frac{(t\omega(z))^{k}}{k!}dz.\label{4.6}\end{equation}
In view of \eqref{4.4} the relation
\begin{equation}\mu(x)-\mu(0)=\frac{1}{2\pi}
\int_{-\infty}^{\infty}\frac{e^{-ixz}-1}{-iz}
\omega(z)dz \label{4.7}\end{equation}holds. The measure $\nu(dx)$ is continuous.Hence formula \eqref{4.7}
implies, that
\begin{equation}
\int_{-\infty}^{\infty}\frac{e^{-ixz}-1}{-iz}
\omega(z)dz\Big{|}_{x=0}=0.\label{4.8}\end{equation} In the same way we can prove the formulas
\begin{equation}
\int_{-\infty}^{\infty}\frac{e^{-ixz}-1}{-iz}
\omega^{k}(z)dz\Big{|}_{x=0}=0,\,k{\geq}1.\label{4.9}\end{equation}
It follows from (4.5) and (4.9) that
\begin{equation}F_{0}(+0,t)-F_{0}(-0,t)=e^{-t\Omega}\lim_{x{\to}+0}\frac{1}{\pi}
\int_{-\infty}^{\infty}\frac{e^{-ixz}-1}{-iz}dz=e^{-t\Omega}.\label{4.10}\end{equation}
The Proposition is proved.\\
We introduce the sequence of functions
\begin{equation}
F_{n+1}(x,t)=\int_{0}^{t}\int_{-\infty}^{\infty}F_{0}(x-\xi,t-\tau)V(\xi)d_{\xi}F_{n}(\xi,\tau)d\tau,\,n{\geq}0,
 \label{4.11}\end{equation} where the function $V(x)$ is defined by relations
 $V(x)=1$ when $x{\notin}\Delta$ and $V(x)=0$ when $x{\in}\Delta.$
 In the right side of  (4.11) we use Stieltjes integration. According to (4.11)the function $F_{1}(x,t)$ is  continuous with respect to $x$, when $x{\ne}0$. The point $x=0$ we shall consider separately.
  \begin{Tm}\label{Theorem 4.3}  Let Conditions 4.1 be fulfilled. If the point $x=0$ belongs to $\Delta$ then
 the functions $F_{n}(x,t),\,n>0$ are continuous with respect to $x$.\end{Tm}
 \emph{Proof.} Using (4.11) we have \begin{equation}F_{1}(+0,t)-F_{1}(-0,t)=\int_{0}^{t}[F_{0}(+0,t-\tau)-F_{0}(-0,t-\tau)]V(0)
 d\tau.\label{4.12}\end{equation}If the point $x=0$ belongs to $\Delta$ then $V(0)=0$.
 Hence,
 $F_{1}(+0,t)-F_{1}(-0,t)=0.$ Thus, the function $F_{1}(x,t)$ is continuous with respect to $x$. Now the assertion of the Theorem follows directly from (4.11).\\
 It follows from (1.1) that
\begin{equation}\mu(z,t)=\mu(z,t-\tau)\mu(z,\tau).\label{4.13}\end{equation}
Due to (4.13) and convolution formula for Stieltjes-Fourier transform ([3], Ch.4)
the relation
 \begin{equation}
F_{0}(x,t)=\int_{-\infty}^{\infty}F_{0}(x-\xi,t-\tau)d_{\xi}F_{0}(\xi,\tau)  \label{4.14}\end{equation} is true.
 Using (4.11) and (4.14) we have
\begin{equation}
0{\leq}F_{n}(x,t){\leq}t^{n}F_{0}(x,t)/n!.\label{4.15} \end{equation}
 Hence the
series \begin{equation}
F(x,t,u)=\sum_{n-0}^{\infty}(-1)^{n}u^{n}F_{n}(x,t)\label{4.16}  \end{equation}
converges. The probabilistic meaning of $F(x,t,u)$ is defined by the
relation (see \cite{Kac1}, Ch.4):
\begin{equation}E\{\mathrm{exp}[-u\int_{0}^{t}V(X_{\tau})d\tau],c_{1}<X_{t}<c_{2}\}=
F(c_2,t,u)-F(c_1,t,u). \label{4.17} \end{equation}
The inequality $V(x){\geq}0$ and
relation (4.17) imply that the function $F(x,t,u)$ monotonically decreases with respect
 to the variable "$u$" and  monotonically increases with respect
 to the variable "$x$". Hence, the following formula
  \begin{equation}
 0{\leq}F(x,t,u){\leq}F(x,t,0)=F_{0}(x,t)\label{4.18}
  \end{equation}
is true. In view of \eqref{4.18} the Laplace transform
\begin{equation}
\Psi(x,s,u)=\int_0^{\infty}e^{-st}F(x,t,u)dt, \quad s>0 \label{4.19}
\end{equation}
%\begin{equation} \Psi(x,s,u)=\int_{0}^{\infty}e^{-st}F(x,t,u)dt,\quad
%s>0. \label(4.19)
%\end{equation}
is correct.
Since the function $F(x,t,u)$ monotonically decreases with respect to $u$,  this is also true for the function $\Psi(x,s,u)$. Hence the limits
\begin{equation} F_{\infty}(x,t)=\lim{F(x,t,u)},\,\Psi_{\infty}(x,s)=\lim{\Psi(x,s,u)},\,
 u{\to}\infty \label{4.20}
 \end{equation}
  exist.
 It follows from \eqref{4.17} that
 \begin{equation}p(t,\Delta)=P(X_{\tau}{\in}\Delta,\,0<\tau<t)=\int_{\Delta}d_{x}F_{\infty}(x,t).
 \label{4.21}\end{equation}
Hence we have
\begin{equation}
 \int_{0}^{\infty}e^{-st}p(t,\Delta)dt=\int_{\Delta}d_{x}\Psi_{\infty}(x,s).  \label{4.22}\end{equation}
2. Relations \eqref{4.22} implies the following assertion
\begin{Pn}\label{Proposition 4.4}Let Conditions 4.1 be fulfilled. If the point $x=0$ belongs to $\Delta$ then
the function $\Psi_{\infty}(x,s)$ for all $s>0$
is monotonically increasing and continuous with respect to $x$ if $x{\ne}0$.\end{Pn}
Using \eqref{4.10} and \eqref{4.19} we have the assertion.
\begin{Pn}\label{Proposition 4.5}Let Conditions 4.1 be fulfilled. If the point $x=0$ belongs to $\Delta$ then
\begin{equation}\Psi(+0,s,u)-\Psi(-0,s,u)=\Psi_{\infty}(+0,s)-\Psi_{\infty}(-0,s)
=\frac{1}{s+\Omega},\,s>0.\label{4.23}\end{equation}\end{Pn}
The behavior of  $\Psi_{\infty}(x,s)$ when $s=0$ we shall consider separately,
using the following Hengartner and Thedorescu result (\cite{HenTheo}, see \cite{Sato}
too):
\begin{Tm}\label{Theorem 4.6} Let $X_t$ be a Levy process.Then for any finite interval K
the estimation
\begin{equation}P(X_t{\in}K)=0(t^{-1/2}) \quad as\quad t{\to}\infty \label{4.24}\end{equation}
is valid.\end{Tm}Hence, we have the assertion (see \cite{Sakh7})
\begin{Tm}\label{Theorem 4.7}Let $X_t$ be a Levy process.Then for any integer ${n>o}$
the estimation
\begin{equation}p(t,\Delta)=0(t^{-n/2}) \quad as\quad t{\to}\infty \label{4.25}\end{equation}
is valid.\end{Tm}

 We need the following partial case of \eqref{4.22}:
\begin{equation}
 \int_{0}^{\infty}p(t,\Delta)dt=\int_{\Delta}d_{x}\Psi_{\infty}(x,0),\label{4.26} \end{equation} According to \eqref{4.25} the integral in the left side of \eqref{4.26} exists.\\
2. \emph{Let us investigate in details the functions $F(x,t,u)$ and $\Psi(x,s,u).$} \\
According to \eqref{4.11} and \eqref{4.16} the function
$F(x,t,u)$ satisfies the equation
\begin{equation}
F(x,t,u)+u\int_{0}^{t}\int_{-\infty}^{\infty}F_{0}(x-\xi,t-\tau)V(\xi)d_{\xi}F(\xi,\tau,u)d\tau=
F_{0}(x,t).\label{4.27}  \end{equation}
Taking from both parts of \eqref{4.27} the Laplace transform, using   \eqref{4.19} and
the convolution property (see \cite{Boch}, Ch.4) we obtain
\begin{equation}
\Psi(x,s,u)+u\int_{-\infty}^{\infty}\Psi_{0}(x-\xi,s)V(\xi)d_{\xi}\Psi(\xi,s,u)=
\Psi_{0}(x,s), \label{4.28} \end{equation}where
\begin{equation}
\Psi_{0}(x,s)=\int_{0}^{\infty}e^{-st}F_{0}(x,t)dt.\label{4.29}  \end{equation}
It follows from \eqref{1.1} and \eqref{4.29}that
\begin{equation}\int_{-\infty}^{\infty}e^{ixp}d_{x}\Psi_{0}(x,s)=\frac{1}{s+\lambda(p)}.
\label{4.30}\end{equation}According to \eqref{4.28} and \eqref{4.29}
 we have \begin{equation}
\int_{-\infty}^{\infty}e^{ixp}[s+\lambda(p)+uV(x)]d_{x}\Psi(x,s,u)=1. \label{4.31}\end{equation}
 Now we introduce the function
\begin{equation}h(p)=\frac{1}{2\pi}\int_{\Delta}e^{-ixp}f(x)dx, \label{4.32}\end{equation}
where the function $f(x)$ belongs to $C_{\Delta}$.\\
By $C_{\Delta}$
 we denote the set of functions $g(x)$ on $L^{2}(\Delta)$ such that
\begin{equation}g(a_{k})=g(b_{k})=g^{\prime}(a_{k})=g^{\prime}(b_{k})=0,\quad 1{\leq}k{\leq}n,\quad g^{\prime\prime}(x){\in}L^{p}(\Delta),\quad p>1.
 \label{4.33}\end{equation} Multiplying both parts  of \eqref{4.31}
by $h(p)$ and integrating them with respect to $p\quad  (-\infty<p<\infty)$ we deduce
the equality
\begin{equation}
\int_{-\infty}^{\infty}\int_{-\infty}^{\infty}e^{ixp}[s+\lambda(p)]h(p)d_{x}\Psi(x,s,u)dp=f(0).
 \label{4.34}\end{equation} We have used the relations \begin{equation}
V(x)f(x)=0,\quad -\infty<x<\infty,\label{4.35}  \end{equation}
\begin{equation}
\frac{1}{2\pi}\mathrm{lim}\int_{-N}^{N}\int_{\Delta}e^{-ixp}f(x)dxdp=f(0),\quad N{\to}\infty.\label{4.36}
  \end{equation}
Since the function $F(x,t,u)$ monotonically decreases with respect
to $"u"$,
 the function $\Psi(x,s,u)$ (see \eqref{4.19}) monotonically decreases with respect
to $"u"$ as well. Hence there exist the limits
 \begin{equation}F_{\infty}(x,s)=\mathrm{\lim}F(x,s,u),\quad
 \Psi_{\infty}(x,s)=\mathrm{\lim}\Psi(x,s,u),\quad u{\to}\infty. \label{4.37} \end{equation}
Using   relations \eqref{1.2} and \eqref{1.6}
 we deduce that
 \begin{equation}\lambda(z)\int_{-\infty}^{\infty}e^{-i{z}\xi}f(\xi)d\xi=
 -\int_{-\infty}^{\infty}e^{-i{z}\xi}[Lf(\xi)]d\xi.\label{4.38}\end{equation}

Relations
\eqref{4.34},\eqref{4.37}  and \eqref{4.38} imply the following assertion.
\begin{Tm}\label{Theorem 4.8} {Let Conditions 4.1 be fulfilled. If the point $x=0$ belongs to $\Delta$
  then the relation  \begin{equation}
  \int_{\Delta}(sI-L_{\Delta})fd_{x}\Psi_{\infty}(x,s)=f(0) \label{4.39} \end{equation} is true.}
 \end{Tm}
\begin{Rk}\label{Remark 4.9}For the Levy processes of type II equality \eqref{4.39} was deduced in the paper
\cite{Sakh7}.\end{Rk}
 Now we shall prove the following assertion.
 \begin{Tm}\label{Theorem 4.11} {Let Conditions 4.1 and inequality \eqref{2.7} be fulfilled. If the point $x=0$ belongs to $\Delta$
  then the  function
  \begin{equation}
  \Psi(x,s)=(I+sB^{\star})^{-1}\Phi(0,x)  ,\label{4.40} \end{equation}
 satisfies relation (4.39).}
  \end{Tm}

\emph{Proof}
In view of (2.8) we have   \begin{equation}
  -BL_{\Delta}f=f,\quad f{\in}C_{\Delta}.\label{4.41} \end{equation}
Relations \eqref{4.40} and \eqref{4.41}
  imply that \begin{equation}
((sI-L_{\Delta})f,\Psi(x,s))_{\Delta}=-((I+sB)L_{\Delta}f,\Psi)_{\Delta}=-(L_{\Delta}f,\Phi(0,x))_{\Delta}. \label{4.42}\end{equation} It is easy to see that
\begin{equation}\Phi(0,x)=B^{\star}\sigma(x),\label{4.43}\end{equation} where $\sigma(x)=-1/2$ when $x<0$ and $\sigma(x)=1/2$ when $x>0$.
 Then according to \eqref{4.40} and \eqref{4.42} relation \eqref{4.39} is
true.\\
 The theorem is proved.

\section{Long time behavior}
\label{sec 5}

1.We apply the following Krein-Rutman theorem (\cite{KR},section 6):\\
\begin{Tm}\label{Theorem 5.1} {If a linear compact operator $T_1$ leaving
invariant a cone $K$,  has a point of the spectrum different from
zero,  then it has a positive eigenvalue $\lambda_{1}$ not less in
modulus than any other eigenvalues $\lambda_{k},\quad (k>1)$.
 To this eigenvalue $\lambda_{1}$ corresponds
at least one eigenvector $g_{1}{\in}K, (T_{1}g_{1}=\lambda_{1}g_{1})$ of
the operator $T_1$ and at least one eigenvector $h_{1}{\in}K^{\star},
(T_{1}^{\star}h_{1}=\lambda_{1}h_{1})$ of the operator $T_{1}^{\star}$.}
\end{Tm}
We remark that in our case the operator $T_1$ has the form\eqref{2.8},the cone $K$ consists of non-negative
continuous real functions $g(x){\in}D_{\Delta}$ and the cone $K^{\star}$ consists of monotonically increasing bounded real functions $h(x){\in}D_{\Delta}^{\star}.$\\
In this section we investigate the asymptotic
behavior $p(t,\Delta)$ when $t{\to}\infty$.\\
2. The spectrum of the operator $B=(1/\Omega)I+T_1$ is situated in the domain $\Re{z}>0.$
The eigenvalue $\mu_1=1/\Omega+\lambda_1$ of the operator $B$ is greater in modulus than any other eigenvalues $\mu_k,\,k>1$ of $B$. We introduce the domain $D_{\ve}$:
\begin{equation}|z-1/{\Omega}|<\ve,\,0<\ve<1/{\Omega}.\label{5.1}\end{equation}
We denote the boundary of the domain  $D_{\ve}$ by $\Gamma_{\ve}$. If $z$ belongs to
$D_{\ve}$ then the relations
\begin{equation}\Re{(1/z)}>c_{\ve}>0 \label{5.2}\end{equation} holds.
We denote
\begin{equation}rank{\lambda_1}=r. \label{5.3}\end{equation}
  Now we formulate the main result of this section.
\begin{Tm}\label{Theorem 5.2}
Let Levy process $X_t$ have  type  $I_c$,\,$0{\in}\Delta$  and let the corresponding
operator $T_1$ satisfy the following conditions:\\
1.Operator $T_1$ is compact in the Banach space $D_{\Delta}.$\\
2. Operator $T_1$ has a point of the spectrum different from zero.\\
Then the  asymptotic equality \begin{equation}
 p(t,\Delta)=e^{-t/\mu_{1}}[q+o(1)],\quad t{\to}+\infty,\quad q{\geq}0
  \label{5.4} \end{equation}
  {is true.}
\end{Tm}
\emph{Proof.}
Using \eqref{4.40}  we obtain   the equality
\begin{equation} p(t,\Delta)=\frac{1}{2\pi}\int_{-\infty}^{\infty}(e^{iyt},\Psi_{\infty}(x,iy))_{\Delta}dy,
\,t>0.\label{5.5}\end{equation} Changing the variable $z=i/y$
 we rewrite \eqref{5.5} in the form
\begin{equation} p(t,\Delta)=\frac{1}{2i\pi}\int_{-i\infty}^{i\infty}(e^{-t/z},(zI-B^{\star})^{-1}
\Phi(0,x))_{\Delta}\frac{dz}{z},
\,t>0.\label{5.6}\end{equation}
As the operator $T_1$ is compact, only a finite number of
eigenvalues $\lambda_{k}, \quad 1<k{\leq}m$ of this operator does
not belong to the domain $D_{\epsilon}$.
 We deduce from formula \eqref{5.6}  the relation
 \begin{equation}
 p(t,\Delta)=
 \sum_{k=1}^{m}\sum_{j=0}^{n_{k}-1}e^{-t/\lambda_{k}}t^{j}c_{k,j}+J, \label{5.7} \end{equation}
where $n_{k}$ is the index of the eigenvalue
$\lambda_{k}$,
\begin{equation}
J=-\frac{1}{2i\pi}\int_{\Gamma_{\ve}}\frac{1}{z}e^{-t/z}(1,(B^{\star} -zI)^{-1}\Phi(0,x))_{\Delta}dz.
 \label{5.8} \end{equation}
 We remind that the \emph{index} of the  eigenvalue
$\lambda_{k}$ is defined as the dimension of the largest Jordan block associated to that eigenvalue.
We note that \begin{equation}n_{1}=1. \label{5.9} \end{equation}
Indeed, if $n_{1}>1$ then there exists such a function $f_{1}$ that
\begin{equation}
Bf_{1}=\lambda_{1}f_{1}+g_{1}. \label{5.10} \end{equation}
In this case the relations
\begin{equation}
(Bf_{1},h_{1})_{\Delta}=\lambda_{1}(f_{1},h_{1})_{\Delta}+(g_{1},h_{1})_{\Delta}
=\lambda_{1}(f_{1},h_{1})_{\Delta}
 \label{5.11} \end{equation}are true. Hence $(g_{1},h_{1})_{\Delta}=0,$ this relation contradicts
 \eqref{5.3}. It proves equality \eqref{5.9}.
 
 Relation (4.43) implies  that
 \begin{equation}
 \Phi(0,x){\in}D^{\star}_{\Delta}. \label{5.12} \end{equation}
Among the numbers $\mu_{k}$ we choose the ones for which
$\mathrm{Re}(1/{\mu}_{k}),\quad (1{\leq}k{\leq}m)$ has the
smallest value $\delta$. Among the obtained numbers we choose
$\mu_{k},\quad (1{\leq}k{\leq}\ell)$ the indexes $n_{k}$ of which
have the largest value $n$. We deduce from \eqref{5.7} and \eqref{5.8} that
\begin{equation}
 p(t,\Delta)=e^{-t{\delta}}t^{n}
 [Q(t)+o(1)],\quad t{\to}\infty. \label{5.13} \end{equation}
 We note that the function
 \begin{equation} Q(t)= \sum_{k=1}^{\ell}e^{it\mathrm{Im}(\mu_{k}^{-1})}c_{k}  \label{5.14} \end{equation}
 is almost periodic (see \cite{LEV}). Hence in view of \eqref{5.13} and the inequality\\
  $p(t,\Delta)>0,\quad (t{\geq}0)$ the following relation
 \begin{equation} Q(t){\geq}0,\quad -\infty<t<\infty  \label{5.15} \end{equation}
 is valid.\\
  First we assume that at least one of the inequalities
 \begin{equation}\delta<{\lambda}_{1}^{-1},\quad n>1 \label{5.16} \end{equation}
 is true. Using  \eqref{5.16} and the inequality
  \begin{equation} \lambda_{1}>|\lambda_{k}|,\quad k=2,3,... \label{5.17} \end{equation}we have
\begin{equation}\mathrm{Im}\mu_{j}^{-1}{\ne}0,\quad 1{\leq}j{\leq}\ell. \label{5.18} \end{equation}
It follows from \eqref{5.14} that
\begin{equation}
c_{j}=\mathrm{lim}\frac{1}{2T}\int_{-T}^{T}Q(t)e^{-it(\mathrm{Im}{\mu}_{j}^{-1})}dt,\quad
T{\to}\infty. \label{5.19} \end{equation}In view of \eqref{5.15}  the relations
\begin{equation}
|c_{j}|{\leq}\mathrm{lim}\frac{1}{2T}\int_{-T}^{T}Q(t)dt=0,\quad
T{\to}\infty, \label{5.20} \end{equation} are valid, i.e. $c_{j}=0,\quad
1{\leq}j{\leq}\ell.$ This means that  relations \eqref{5.14} are not true.
Hence the equalities
 \begin{equation}\delta={\lambda}_{1}^{-1},\quad n=1 \label{5.21} \end{equation}
 are true. From  \eqref{5.21} we get the asymptotic equality
 \begin{equation}
 p(t,\Delta)=e^{-t/\lambda_{1}}[q+o(1)]\quad t{\to}\infty, \, q{\geq}0,\label{5.22} \end{equation}
where
 \begin{equation} q=\sum_{k=1}^{r}g_{k}(0)\int_{\Delta}d{h_{k}(x)}.
 \label{5.23} \end{equation}
  Here $g_{k}(x)$ are the eigenfunctions of the operator $B$ corresponding to
  the eigenvalues $\mu_{1}$, and $h_{k}(x)$ are the eigenfunctions of the operator $B^{\star}$
   corresponding to
  the eigenvalues $\mu_{1}$ The following conditions are fulfilled
  \begin{equation}
 (g_{k},h_{k})=\int_{\Delta}g_{k}(x)dh_{k}(x)=1, \label{5.24} \end{equation}
\begin{equation}
 (g_{k},h_{\ell})=\int_{\Delta}g_{k}(x)dh_{\ell}(x)=0,\quad k{\ne}\ell. \label{5.25}\end{equation} T he theorem is proved.

\section{Example}
1.Let us consider the example, when
\begin{equation}\nu^{\prime}(x)=p^{2}e^{-p|x|},\,p>0,\,\Delta=[0,\omega],
\,A=0,\,\gamma=0.\label{6.1}\end{equation}
Using (1.26) we have
\begin{equation}L_{\Delta}f=-2p[f(x)-(p/2)\int_{0}^{\omega}e^{-p|x-y|}f(y)dy].\label{6.2}\end{equation}
Condition (2.6) is fulfilled. Hence the operator $L_{\Delta}^{-1}$ has form
(2.7), where the operator $T_1$ is defined by the relation
\begin{equation}T_{1}f=\int_{0}^{\omega}\gamma(x,t)f(t)dt.\label{6.3}\end{equation}
It follows from (2.7) and (6.2), that
\begin{equation}\gamma(x,t)-(p/2)e^{-p|x-t|}-(p/2)\int_{0}^{\omega}e^{-p|x-y|}\gamma(y,t)dy=0.
\label{6.4}\end{equation}According to (6.4) we have
\begin{equation}\frac{\partial^2{\gamma}}{\partial{x^{2}}}=0,\,x{\ne}t.\label{6.5}\end{equation}
From \eqref{6.5} we obtain that
\begin{equation}\gamma(x,t)=c_1(t)+c_2(t)x,\, x>t,\label{6.6}\end{equation}
\begin{equation}\gamma(x,t)=c_3(t)+c_4(t)x,\, x<t.\label{6.7}\end{equation}
With the help of \eqref{6.6}, \eqref{6.7} and equality $\gamma(x,t)=\gamma(t,x)$ we deduce that
\begin{equation}\gamma(x,t)=(\alpha_{1}+\alpha_{2}t)+(\beta_{1}+\beta_{2}t)x,\,t<x.
\label{6.8}\end{equation}
Using \eqref{6.6} and relations
\begin{equation}
\gamma(x,0)-(p/2)e^{-px}-(p/2)\int_{0}^{\omega}e^{-p|x-y|}\gamma(y,0)dy=0,
\label{6.9}\end{equation}
\begin{equation}\frac{\partial\gamma(x,t)}{\partial{t}}\Big|_{t=0}-(p^{2}/2)e^{-px}-(p/2)\int_{0}^{\omega}e^{-p|x-y|}
\frac{\partial\gamma(y,t)}{\partial{t}}\Big|_{t=0}dy=0.
\label{6.10}\end{equation} we obtain
\begin{equation}\alpha_1=\frac{p(1+p\omega)}{2+p\omega},\,\beta_1=-\frac{p^2}{2+p\omega},\,
\alpha_2=p\alpha_1,\,\beta_2=p\beta_1.\label{6.11}\end{equation}
Thus, the corresponding operator $T_1$ has the form \eqref{6.3}, where the kernel $\gamma(x,t)$ is defined by the relations \eqref{6.8} and \eqref{6.11}.\\
2. Let us find the eigenvalues and eigenfunctions of the operator $T_1$, when \eqref{6.1} is valid. The eigenfunctions $f(x)$ and eigenvalues $\lambda$ satisfy the relations
\begin{equation}\int_{0}^{\omega}\gamma(x,t)f(t)dt=\lambda{f(x)},\quad
\int_{0}^{\omega}\frac{\partial}{\partial{t}}\gamma(x,t)f(t)dt=\lambda{f^{\prime}(x)}.
\label{6.12}\end{equation}
It follows from \eqref{6.12} that
\begin{equation}\lambda{f^{\prime\prime}(x)}=-p^{2}f.\label{6.13}\end{equation}
Hence,
\begin{equation}f(x,\lambda)=c_{1}(\lambda)\sin(xp/\sqrt{\lambda})+c_{2}(\lambda)\cos(xp/\sqrt{\lambda}).
\label{6.14}\end{equation}
According to \eqref{6.12} and \eqref{6.14}
we have
\begin{equation}
c_{1}(\lambda)a(\lambda)+c_{2}(\lambda)b(\lambda)=
\lambda{c_{2}(\lambda)},
\label{6.15}\end{equation}
\begin{equation}
c_{1}(\lambda)a(\lambda)+c_{2}(\lambda)b(\lambda)=
\sqrt{\lambda}{c_{1}(\lambda)},\label{6.16}
\end{equation}
where
\begin{equation}a(\lambda)=\frac{1}{2+p\omega}[-\sqrt{\lambda}\cos(p\omega/\sqrt{\lambda}) -{\lambda}\sin(p\omega/\sqrt{\lambda})+(1+p\omega)\sqrt{\lambda}],
\label{6.17}\end{equation}
\begin{equation}b(\lambda)=\frac{1}{2+p\omega}[\sqrt{\lambda}\sin(p\omega/\sqrt{\lambda}) -{\lambda}\cos(p\omega/\sqrt{\lambda})+\lambda].
\label{6.18}\end{equation}
The eigenvalues of the operator $B$ are defined by the relations \eqref{6.15} and \eqref{6.16},i.e.
\begin{equation}\sqrt{\lambda}a(\lambda)+b(\lambda)=\lambda.\label{6.19}\end{equation}
Equalities \eqref{6.17}-\eqref{6.19} imply
\begin{equation}\tan(p\omega/{\sqrt{\lambda}})=\frac{2\sqrt{\lambda}}{1-\lambda}.
\label{6.20}\end{equation}The following table gives the maximal positive roots of the equation \eqref{6.20}.
\begin{equation}\left(
  \begin{array}{cccccc}
    p\omega & \pi/4 & \pi/3 & \pi/2 & 2\pi/3 & \pi \\
    \lambda & 0 .445 & 0.617 & 0.162 & 1.433 & 2.454\\
  \end{array}
\right)\label{6.21}\end{equation}
\section{Discrete Levy measure}
We consider the case when Levy measure $\nu(u)$ is discrete, where
the  gap points we denote by $\nu_{k}\,(1{\leq}k{\leq}n)$ and  the corresponding
 gaps we denote by $\sigma_{k}>0$.  Let condition \eqref{1.7} and \eqref{1.9} be fulfilled. In this case formula \eqref{1.26}
 takes the form
 \begin{equation}L_{\Delta}f=-{\Omega}f(x)+\sum_{k=1}^{n}f(\nu_{k}+x)\sigma_{k},\,
 \Omega{\geq}\sum_{k=1}^{n}\sigma_{k}\label{7.1}\end{equation}
 and the operator $T$ is defined by the relation
 \begin{equation}Tf=\sum_{k=1}^{n}f(\nu_{k}+x)\sigma_{k},\,x{\in}\Delta.
 \label{7.2}\end{equation}
 We note that $f(x)=0 \quad if \quad
 x{\notin}\Delta$. Relation \eqref{7.2} implies the  assertion.\\
 \begin{Pn}\label{Proposition 7.1} If  $\nu_{k}>0,\,(1{\leq}k{\leq}n)$ then  the operator $T$ has the following properties:\\
 1.There exists such integer number m that
  the equality
 \begin{equation}T^{m}=0\label{7.3}\end{equation}
 is valid.\\
 2.The operator $T$ maps the bounded functions in the bounded functions.\\
 3.The operator $T$ maps the non-negative  functions in the non-negative functions.\\
 4.The operator $T_1$ in equality \eqref{2.8} has the form
 \begin{equation}T_{1}=T/\Omega+(T/\Omega)^{2}+...+(T/\Omega)^{m-1}.\label{7.4}\end{equation}
 \end{Pn}
 {\bf Acknoledgements }  The author is very grateful to I. Tydniouk for his calculation of the table \eqref{6.21}.

%%%%%%%
\end{document}